\documentclass[11pt]{article}
\usepackage{amssymb,latexsym,amsfonts,verbatim,amscd}
\usepackage[all]{xy}
\newtheorem{Def}{Definition}[section]
\newtheorem{Thm}[Def]{Theorem}
\newtheorem{Lem}[Def]{Lemma}

\newtheorem{Cor}[Def]{Corollary}

\newtheorem{Que}[Def]{Question}

\newtheorem{Fac}[Def]{Fact}

\font\nat msbm10 scaled\magstephalf

\def\P{\hbox{\nat\char80}}

\def\telos{\hfill$\dashv$}
 \font\goth eufm10

\def\int{\hbox{\rm int}}

\begin{document}
\sloppy

\title{Typicality \`{a} la Russell in set theory}
\author{Athanassios Tzouvaras}

\date{}
\maketitle

\begin{center}
Aristotle University of Thessaloniki \\
Department  of Mathematics  \\
541 24 Thessaloniki, Greece. \\
e-mail: \verb"tzouvara@math.auth.gr"
\end{center}

\begin{abstract}
We adjust   the notion  of typicality originated with  Russell, which was   introduced and studied in a previous paper for general first-order structures,  to  make it expressible in the language of set theory. The adopted definition of the class ${\rm NT}$  of  nontypical sets comes out as  a natural strengthening of Russell's initial definition, which employs properties of small (minority) extensions, when the latter are restricted to  the various levels $V_\zeta$ of $V$. This strengthening leads to defining ${\rm NT}$ as the class of sets  that belong to some countable ordinal definable set.  It follows that  ${\rm OD}\subseteq {\rm NT}$ and hence ${\rm HOD}\subseteq {\rm HNT}$. It is proved that the class  ${\rm HNT}$ of hereditarily nontypical sets  is an inner model of ${\rm ZF}$. Moreover the (relative) consistency of $V\neq {\rm NT}$ is established, by showing that in many forcing extensions $M[G]$ the generic set $G$ is a typical element of  $M[G]$, a fact which is fully in accord  with the intuitive meaning of typicality. In particular it is consistent that there exist continuum many typical reals.
In addition it follows from a result of Kanovei and Lyubetsky that  ${\rm HOD}\neq {\rm HNT}$ is also relatively  consistent. In particular it is consistent that ${\cal P}(\omega)\cap {\rm OD}\subsetneq{\cal P}(\omega)\cap {\rm NT}$.  However many questions remain  open, among them  the consistency of  ${\rm HOD}\neq {\rm HNT}\neq V$, ${\rm HOD}={\rm HNT}\neq V$ and  ${\rm HOD}\neq {\rm HNT}= V$.

\end{abstract}

\vskip 0.2in


\vskip 0.2in


\section{Introduction}
Inspired by B. Russell's definition of ``typical Englishman'' in \cite{Ru95}, we introduced in \cite{Tz20} a rigorous definition of {\em typical property} over a first-order structure $\mbox{\goth M}=(M,\ldots)$ (expressed in the language $L$ of $\mbox{\goth M}$), and of {\em typical element} of $\mbox{\goth M}$. The definition of typical element exhibits  clear similarities  with the well-known definition of Martin-L\"{o}f random real and its variants  so it could be conceived as  an alternative notion of randomness. Moreover, allowing parameters in the formulas involved leads to a notion of {\em relative typicality} ``object $a$ is typical with respect to objects $\overline{b}$'', or ``$\overline{b}$-typical'', denoted $\textsf{Tp}(a,\overline{b})$. In this form, typicality
exhibits also clear similarities with van Lambalgen's relative randomness notion $R(x,\overline{y})$, and $\textsf{Tp}(a,\overline{b})$ was shown to satisfy most of the randomness axioms proposed in  \cite{La90} and other papers of the same author.  Let us first recall some basic definitions from \cite{Tz20}.

Russell's inspiring  definition of typical Englishman \cite[p. 89]{Ru95} is as follows:
\begin{quote}
{\em A typical Englishman is one who possesses  all the properties possessed by a majority of Englishmen.}
\end{quote}
Without the distinction between object language and metalanguage, the definition is obviously circular, so practically useless, but it becomes natural and sound as soon as we make the aforementioned distinction. This is easily made by the help of elementary tools of model theory.

Throughout our background theory is ${\rm ZFC}$. Let $L$ be  a first-order language, $\mbox{\goth M}=(M,\ldots)$  an infinite $L$-structure and $A\subseteq M$. $L(A)$ denotes $L$ augmented with parameters from $A$. By some abuse of language we refer also to $L(A)$ as the set of formulas of $L(A)$.   By a property of $L(A)$ we mean  a formula $\phi(x)\in L(A)$ with  one free variable. A  {\em majority  set} is any $X\subseteq M$ that contains strictly more elements than its complement, i.e.,  $|X|>|M\backslash X|$. (Accordingly $M\backslash X$ is then a {\em minority set}).  Note that   $|X|>|M\backslash X|$ is equivalent to  $|M\backslash X|<|M|$ (which in particular  implies that $|X|=|M|$). Let
$$\textsf{mj}(M)=\{X\subseteq M:|X|>|M\backslash X|\}=\{X:|M\backslash X|<|M|\}$$
be the class of majority subsets of $M$.
$\textsf{mj}(M)$ is a filter on $M$ that  extends the Fr\'{e}chet filter of cofinite subsets of $M$.

Given an $L$-structure  $\mbox{\goth M}$ and  a property $\phi(x,\overline{y})$ of $L$ and parameters $\overline{b}\in M$, $\phi({\mbox{\goth M}},\overline{b})$  denotes the extension of $\phi(x,\overline{b})$ in $\mbox{\goth M}$, i.e., $$\phi({\mbox{\goth M}},\overline{b})=\{a\in M:\mbox{\goth M}\models\phi(a,\overline{b})\}.$$

\begin{Def} \label{D:oldtypical}
{\em    A property $\phi(x,\overline{y})$ of $L$ is said to be} $A$-typical over $\mbox{\goth M}$, {\em  for some $A\subseteq M$, if for every $\overline{b}\in A$, $\phi({\mbox{\goth M}},\overline{b})\in \textsf{mj}(M)$. In particular  $\phi(x)$ is} typical {\em  if it is $\emptyset$-typical.  An element $a\in M$ is said to be} $A$-typical {\em if it satisfies every $A$-typical property over $\mbox{\goth M}$.}
\end{Def}

If $A$ is finite we write it in the form   $\overline{a}=\langle a_1,\ldots,a_n\rangle$, and say  $\overline{a}$-typical instead of $A$-typical. Given a structure $\mbox{\goth M}$, the notation $\textsf{Tp}(a,\overline{b})$  means ``$a$ is $\overline{b}$-typical'' over $\mbox{\goth M}$. $\textsf{Tp}(a)$ means  $a$ is ``$\emptyset$-typical'', or just ``typical''.

The following existence results for typical elements were shown in \cite{Tz20}.

\begin{Thm} \label{T:suffices}
If $\mbox{\goth M}$ is $\kappa$-saturated, for some $\kappa\geq \aleph_0$, then for every $A\subseteq M$ such that $|A|<\kappa$, $M$ contains $A$-typical elements.
\end{Thm}

\begin{Thm} \label{T:typ}
Let $\mbox{\goth M}$ be an $L$-structure, for a countable $L$, and let $A\subseteq M$ be a set of parameters  such that $cf(|M|)>\max(\aleph_0,|A|)$.  Then $\mbox{\goth M}$ contains $|M|$ $A$-typical elements, while only $<|M|$ non-$A$-typical ones.
\end{Thm}

Applying Theorem \ref{T:typ} to the standard model of full second-order arithmetic  ${\cal R}=(\omega,{\cal P}(\omega),+,\cdot,<,\in,0,1)$ we have the following immediate consequence.

\begin{Thm} \label{T:typicalreal}
For every $A\subseteq {\cal P}(\omega)$ such that $|A|\leq \aleph_0$, there exist $2^{\aleph_0}$ $A$-typical reals, while only $<2^{\aleph_0}$ non-$A$-typical ones. More precisely:  For every finite (or even countable) tuple $\overline{b}$ of reals, $|\{x:\textsf{Tp}(x,\overline{b})\}|=2^{\aleph_0}$, while  $|\{x:\neg \textsf{Tp}(x,\overline{b})\}|<2^{\aleph_0}$.
\end{Thm}

\section{Typicality in ZFC}
The definitions and facts given in the Introduction refer to  what might be called ``external typicality'' with respect to a given first-order structure $\mbox{\goth M}$, because the cardinality measure used to tell which properties are typical and which are not, is external with respect to $\mbox{\goth M}$.  In  \cite{Tz20} (Remark 2.15) we already observed that  such a concept is suitable for all first-order structures except for models of set theory, because every such model  possesses its own  internal cardinality, which is of course   preferable  over the external one, and we stressed ``the challenge to find a notion of typicality suitable for  models of set theory.'' It is the purpose of this paper to examine  ways for implementing  Russell's notion of typicality naturally in any model of ${\rm ZFC}$. This amounts to finding  a typicality notion expressible in the language of set theory.  The basic problem has to do with the implementation of Russell's   majority/minority criterion concerning the size of extensions of properties in the universe $V$ of ${\rm ZFC}$.

We have argued in \cite{Tz20} that typical objects  behave  much like {\em random} entities, which is equivalent to saying that nontypical objects are  expected to  be ``special'' entities, of the kind of objects by which we build  ``familiar'' structures. For instance all definable elements of a structure are nontypical according to the definitions of the previous section (referring to external typicality), and clearly the same must be true for any  reasonable notion of typicality in  ${\rm ZFC}$. That is to say, every set  definable without parameters in the universe $V$ should be nontypical. To simplify notation, henceforth we shall identify a structure $\mbox{\goth M}$ with its domain $M$. Given a set or class structure $M$,  we denote by ${\rm Df}(M)$ the set or class  of first-order definable (without parameters) elements of $M$ (i.e., $a\in {\rm Df}(M)$ if there is a formula $\phi(x)$ such that $a$ is the unique element of $M$ that satisfies $\phi(x)$). In particular,  ${\rm Df}(V)$ is the class of all sets definable without parameters in the universe $(V,\in)$. By Tarski's undefinability of truth,  ${\rm Df}(M)$ is not a definable subclass of $M$ (in particular, the fact that  a set belongs to ${\rm Df}(V)$ is not expressible by a formula of set theory). On the other hand it is of some interest to note  (although we are not going to use this fact below) that ${\rm Df}(V)$ is an elementary subclass of the inner model HOD of hereditarily definable sets of ${\rm ZFC}$, i.e., ${\rm Df}(V)\preccurlyeq{\rm HOD}$ (see e.g. \cite[Thm. 4]{HLR13}).

Let ${\rm TP}$, ${\rm NT}$  denote the complementary classes of {\em typical} and {\em nontypical} sets of $V$, respectively, the definitions of which are being sought.  It follows from the foregoing discussion that  the definition of ${\rm NT}$ should satisfy the following requirement:
\begin{equation} \label{E:first}
{\rm Df}(V)\subseteq {\rm NT}.
\end{equation}
Further, as already mentioned above,  the definition of  ${\rm NT}$ should  be  expressible in the language of set theory. Recall from the previous section that given a first structure $M$ and a formula $\phi(x)$ (or $\phi(x,\overline{b})$, with  $\overline{b}\in M$), $\phi(M)$ (resp. $\phi(M,\overline{b})$) denotes  the extension of  $\phi(x)$ in $M$.
The same notation will be  used in set theory even when $M$ is a proper class. So  $\phi(V)$ is the class $\{x:(V,\in)\models\phi(x)\}$, which is commonly written   $\{x:\phi(x)\}$. Thus a second requirement for the definition of  ${\rm NT}$ is that there is a first-order property $\sigma(x)$ (without parameters) such that:
\begin{equation} \label{E:second}
{\rm NT}=\sigma(V).
\end{equation}
Apart from requirements (\ref{E:first}) and (\ref{E:second}), the definition  should be nontrivial, e.g.  ${\rm ZFC}$ should not prove that there are no typical elements (i.e., $V={\rm NT}$), or that ${\rm NT}$ coincides with some already known subclass of $V$.

Below we  examine three  options for the definition of ${\rm NT}$.  All  of them  are  attempts to implement the Russell' majority/minority criterion for extensions of properties of a set $a$. The first option  leads to $V={\rm NT}$ and thus, inevitably,  is rejected as trivial. The second one lies, in a sense, to the opposite end as it makes  $V\neq {\rm NT}$  provable in ${\rm ZFC}$. A consequence of this is that  ${\rm Df}(V)\not\subseteq {\rm NT}$, i.e., the refutation of  condition (\ref{E:first}), is consistent. This is a sufficient reason to reject this  option as well.  The third option leads to  two sub-options, one of which, happily, offers a definition  for ${\rm NT}$ according to which its subclass ${\rm HNT}$ (of hereditarily nontypical sets) turns out to be a {\em new}  inner model of ${\rm ZF}$.

Since dealing with extensions $\phi(V)$ over the entire $V$ leads to notions non-expressible in the language of set theory,  we must be confined to  extensions of $\phi$'s in the segments $V_\zeta=\{x:rank(x)<\zeta\}$, $\zeta\in Ord$, of $V$, i.e., to sets $\phi(V_\zeta)=\{x:(V_\zeta,\in)\models\phi(x)\}$. (We  shall write simply $V_\zeta\models\phi$ instead of $(V_\zeta,\in)\models\phi$.)  Working in $(V_\zeta,\in)$, we often use ordinal parameters from $V_\zeta$, i.e., tuples  $\overline{\theta}<\zeta$, but, as will be explained later, the collection of sets defined {\em over all} $(V_\zeta,\in)$ using ordinal parameters is no different than the collection of sets defined in the same structures with no parameters at all.  So we can safely dispense, at least for the moment, with ordinal parameters.

Before going on let us recall the formal definition of Russell's typicality and  fix some notation. Given a structure $M$,  let
$${\rm TP}(M)=\{a\in M:\textsf{Tp}(a)\}, \quad {\rm NT}(M)=\{a\in M:\neg\textsf{Tp}(a)\}$$
be the classes of typical and nontypical elements of $M$, respectively. It follows   from Definition  \ref{D:oldtypical} that
\begin{equation} \label{E:explicit2}
{\rm NT}(M)=\{a\in M:(\exists \phi)(a\in \phi(M) \ \wedge \ |\phi(M)|<|M|)\}.
\end{equation}
We come to the three options for the definition of ${\rm NT}(V):={\rm NT}$.

\vskip 0.1in

\textbf{Option 1.} Given that the definition of Russell's typicality does not apply directly to the entire $V$, a natural alternative is to define it with respect to the set-approximations $V_\zeta$, $\zeta\in Ord$, of $V$. Namely, we say that a set $a$ is nontypical if and only if it is so  {\em with respect to  some $V_\zeta$}  to which it belongs (that is, for some  $\zeta>rank(a)$), i.e., we set
\begin{equation} \label{E:option1}
{\rm NT}:=\bigcup_{\zeta>\omega}{\rm NT}(V_\zeta),
\end{equation}
where ${\rm NT}(V_\zeta)$ is defined as in  (\ref{E:explicit2}), that is
$${\rm NT}(V_\zeta)=\{a\in V_\zeta:(\exists \phi)(a\in \phi(V_\zeta)\wedge |\phi(V_\zeta)|<|V_\zeta|)\}.$$
It is easy to see that  definition (\ref{E:option1}) satisfies both conditions (\ref {E:first}) and (\ref {E:second}), but, unfortunately, is trivial in the sense that it implies that  every set is nontypical,  i.e., ${\rm NT}=V$. Indeed, take any  $a\in V$ and pick some  $\zeta$ such that  $a\in V_\zeta$.  Then  clearly $\zeta$ is definable in $V_{\zeta+1}$ (as its  greatest ordinal), i.e., $\zeta\in {\rm Df}(V_{\zeta+1})$, and thus $V_\zeta\in {\rm Df}(V_{\zeta+1})$ (as the set of elements of rank $<\zeta$). Thus   there is $\phi(x)$ such that $V_\zeta=\phi(V_{\zeta+1})$, hence  $a\in \phi(V_{\zeta+1})$. Moreover $|\phi(V_{\zeta+1})|=|V_\zeta|<|V_{\zeta+1}|$. Therefore  $a\in{\rm NT}(V_{\zeta+1})$, and hence $a\in{\rm NT}$. It follows that no typical sets exist, or ${\rm TP}=\emptyset$.

This fact would  be normal and acceptable if it concerned only {\em some  model} of ${\rm ZFC}$. But being a theorem of ${\rm ZFC}$ trivializes ${\rm NT}$, indicating that its definition is too generous. So next we shall  attempt  to restrict it somehow, and widen, accordingly, the class ${\rm TP}$.

\vskip 0.1in

\textbf{Option 2.} We saw above that if $rank(a)=\zeta$, i.e., $a$ first occurs in $V_{\zeta+1}$,  then already $a\in {\rm NT}(V_{\zeta+2})$ and so $a\in {\rm NT}$ according to the definition in Option 1. Thus  a possible way to restrict ${\rm NT}$ and avoid ${\rm NT}=V$ might be to  decide the typicality of a set $a$, not in {\em any level} $V_\zeta$ with $\zeta>rank(a)$, but rather at the first level of its occurrence, i.e., at $V_{rank(a)+1}$. That would lead  to replacing definition  (\ref{E:option1}) of ${\rm NT}$ with the following:
\begin{equation} \label{E:option2}
{\rm NT}_0:=\{a:a\in {\rm NT}(V_{rank(a)+1})\}.
\end{equation}
Then for every  $a$ with $rank(a)=\zeta$,  $a\in {\rm NT}_0$ if and only if there is $\phi(x)$ such that $a\in \phi(V_{\zeta+1})$ and $|\phi(V_{\zeta+1})|<|V_{\zeta+1}|$.

It turns out, however, that definition (\ref{E:option2}) is  now too  restrictive for the class of nontypical sets (and accordingly too generous for the class of typical sets) and results in the (consistency of the) failure of condition (\ref{E:first}). This will be a consequence  of  the following result. (Note that ${\rm NT}(V_\zeta)$ and ${\rm TP}(V_\zeta)$ keep the meaning specified in definition (\ref{E:explicit2}).)

\begin{Lem} \label{L:generous}
{\rm (${\rm ZFC}$)} (i) For every $\zeta\geq\omega$, $|{\rm NT}(V_{\zeta+1})|<|V_{\zeta+1}|$, while $|{\rm TP}(V_{\zeta+1})|=|V_{\zeta+1}|$.

(ii) Moreover, for every $\zeta\geq\omega$, $|{\rm TP}(V_{\zeta+1})\cap(V_{\zeta+1}\backslash V_\zeta)|=|V_{\zeta+1}|$.
\end{Lem}

{\em Proof.} (i) The proof is easy and quite similar to the proof of Theorem \ref{T:typ} above (see \cite{Tz20}, Theorem 2.9). It follows simply from two facts: first, there are countably many properties $\phi(x)$ of $L$ without parameters, and second, the cofinality of $|V_{\zeta+1}|=2^{|V_{\zeta}|}$ is always uncountable. Just observe that  by definition
$${\rm NT}(V_{\zeta+1})=\bigcup\{\phi(V_{\zeta+1}):\phi\in L, \ \mbox{and} \ |\phi(V_{\zeta+1})|<|V_{\zeta+1})|\}.$$
The sets $\phi(V_{\zeta+1})$ on the right-hand side are countably many and  each with cardinality $<|V_{\zeta+1}|$, so, since ${\rm cf}(|V_{\zeta+1}|)>\aleph_0$, $|{\rm NT}(V_{\zeta+1})|<|V_{\zeta+1}|$. Consequently $|{\rm TP}(V_{\zeta+1})|=|V_{\zeta+1}\backslash {\rm NT}(V_{\zeta+1})|=|V_{\zeta+1}|$.

(ii) We have $${\rm TP}(V_{\zeta+1})=({\rm TP}(V_{\zeta+1})\cap V_{\zeta})\cup ({\rm TP}(V_{\zeta+1})\cap (V_{\zeta+1}\backslash V_{\zeta})).$$
Since $|{\rm TP}(V_{\zeta+1})\cap V_{\zeta}|\leq |V_\zeta|$ while, by (i), $|{\rm TP}(V_{\zeta+1})|=|V_{\zeta+1}|$, the claim follows. \telos

\begin{Cor} \label{C:toomany}
${\rm ZFC}$ proves that $V\neq {\rm NT}_0$. Specifically, for every successor level $V_{\zeta+1}$, $|V_{\zeta+1}\cap {\rm TP}_0|=|V_{\zeta+1}|$ (where ${\rm TP}_0=V\backslash {\rm NT}_0$).
\end{Cor}

{\em Proof.} The elements of $V_{\zeta+1}\backslash V_\zeta$ are exactly those of rank $\zeta$, so by the definition of ${\rm TP}_0$, clearly $V_{\zeta+1}\cap {\rm TP}_0=(V_{\zeta+1}\backslash V_{\zeta})\cap {\rm TP}(V_{\zeta+1})$. Thus the claim follows from Lemma \ref{L:generous} (ii).  \telos

\begin{Cor} \label{C:onefails}
${\rm Df}(V)\not\subseteq {\rm NT}_0$  is consistent with  ${\rm ZFC}$. Therefore  condition (\ref{E:first}) is not provable.
\end{Cor}

{\em Proof.} It is known that if ${\rm ZFC}$ is consistent, then it has models all  elements of which are definable, i.e., ${\rm Df}(V)=V$ holds there (see \cite{HLR13} where the existence of such models,  called pointwise definable, is proved). Let $M\models{\rm ZFC}$ be such a model. Then, $M\models V\neq {\rm NT}_0$,  since Corollary \ref{C:toomany} is provable  in ${\rm ZFC}$, or equivalently $M\models V\not\subseteq {\rm NT}_0$. But since  also in $M$ ${\rm Df}(V)=V$, it follows that in this model  ${\rm Df}(V)\not\subseteq {\rm NT}_0$. \telos

\vskip 0.2in

To sum up: replacing definition  (\ref{E:option1}) with definition (\ref{E:option2}) leads to the opposite end of the spectrum, narrowing too much the class of nontypical sets and leading to the non-provability of condition (\ref{E:first}). Therefore some intermediate solution is needed, and this is examined  in the next option.

\vskip 0.1in

\textbf{Option 3.} Our next attempt to restrict ${\rm NT}$ of Option 1 is by requiring a nontypical element of $V_\zeta$ to be caught in the extension $\phi(V_\zeta)$ of some property $\phi$ with  cardinality not just  $<|V_\zeta|$ but rather $<\kappa$, for some $\kappa$ with  $\aleph_0\leq \kappa\leq|V_\zeta|$.

To be a little bit more general, given a structure $M$ and a cardinal $\kappa$ with $\aleph_0\leq\kappa\leq |M|$,  the definition  (\ref{E:explicit2}) of   ${\rm NT}(M)$ can be refined as follows:
\begin{equation} \label{E:kappa2}
{\rm NT}_\kappa(M)=\{a\in M:(\exists \phi)(a\in \phi(M)\ \wedge \ |\phi(M)|<\kappa)\}.
\end{equation}
(In this notation  ${\rm NT}(M)={\rm NT}_{|M|}(M)$.)

\begin{Fac} \label{F:widen}
For any $M$ and any $\aleph_0\leq \kappa<\lambda\leq |M|$, ${\rm NT}_\kappa(M)\subseteq {\rm NT}_\lambda(M)$.   In particular, ${\rm NT}_{\aleph_0}(M)$ is the set of algebraic elements of $M$.
\end{Fac}

{\em Proof.} The first claim follows immediately from  definition  (\ref{E:kappa2}). Concerning the second claim, recall that an element $a\in M$ is {\em algebraic} in $M$ if there is a formula $\phi(x)$ in the language of $M$ without parameters such that $\phi(M)$ is finite and $a\in \phi(M)$. (In particular every element of $M$ which is definable without parameters is algebraic.) Therefore $a$ is algebraic if and only if for some $\phi$, $a\in \phi(M)$ and $|\phi(M)|<\aleph_0$, i.e., if and only if $a\in {\rm NT}_{\aleph_0}(M)$. \telos

\vskip 0.2in

It follows that the class of nontypical elements of $M$ is {\em minimized}, becoming  identical to the class of algebraic elements,  if we set ${\rm NT}(M):={\rm NT}_{\aleph_0}(M)$  (when accordingly the class of typical elements is maximized).
For $M=(V_\zeta,\in)$ and for $\kappa$ such that $\aleph_0\leq \kappa\leq |V_\zeta|$, we have in particular:
\begin{equation} \label{E:parkappa2}
{\rm NT}_\kappa(V_\zeta)=\{a\in V_\zeta:(\exists \phi)(a\in \phi(V_\zeta)\ \wedge \ |\phi(V_\zeta)|<\kappa)\}.
\end{equation}
Applying Fact \ref{F:widen} to the structures $V_\zeta$, we get:

\begin{Fac} \label{F:wideapp}
For any $\zeta>\omega$ and any $\aleph_0\leq \kappa<\lambda\leq |V_\zeta|$,  ${\rm NT}_\kappa(V_\zeta)\subseteq {\rm NT}_\lambda(V_\zeta)$). In particular, ${\rm NT}_{\aleph_0}(V_\zeta)$ is the set of algebraic elements of the structure $(V_\zeta,\in)$.
\end{Fac}
The definition of ${\rm NT}_\kappa(V_\zeta)$ generalizes the definition of ${\rm NT}$ locally but it is not clear if it applies  to the class ${\rm NT}$ itself. Recall that for any ordinal $\alpha\geq 0$, $|V_{\omega+\alpha}|=\beth_\alpha$, where $\beth_0=\aleph_0$, $\beth_{\alpha+1}=2^{\beth_\alpha}$ and $\beth_\alpha=\sup\{\beth_\beta:\beta<\alpha\}$, for a limit $\alpha$. (Under ${\rm GCH}$, $\beth_\alpha=\aleph_\alpha$, for all $\alpha$.)

From the perspective of this paper the  dichotomy of sets into typical and nontypical concerns exclusively infinite sets, or finite sets containing infinite sets etc, while the  elements of  $V_\omega$, being  hereditarily finite, are definable and hence nontypical with respect to any definition of ${\rm NT}$. So the ordinal $\zeta$ in (\ref{E:parkappa2}) must   range beyond $\omega$, i.e., $\zeta>\omega$. Consequently  for each infinite cardinal $\kappa$ it is natural to set
$${\rm NT}_\kappa=\bigcup_{\zeta>\omega}{\rm NT}_\kappa(V_\zeta)$$ as a generalization of (\ref{E:option1}). The question is: what should be the range of $\kappa$? If $|V_\zeta|<\kappa$, for some $\zeta$, then by  (\ref{E:parkappa2}) clearly    $V_\zeta={\rm NT}_\kappa(V_\zeta)$, and hence $V_\zeta\subseteq {\rm NT}_\kappa$.  But given that $\zeta>\omega$, this  seems quite unnatural. For in that case  a whole  segment of $V$ that strictly extends $V_\omega$ should consist exclusively of nontypical sets.  In particular we would have  ${\cal P}(\omega)\subseteq V_{\omega+1}\subseteq {\rm NT}_\kappa$, i.e., all reals would  be nontypical, contrary to the fact that typicality was introduced in \cite{Tz20}  (rather successfully) as a parallel notion of randomness, with a large amount of reals to be proved  typical (with respect to the ``external'' notion of typicality employed there).

The only way to avoid this unnatural  situation is to restrict   ourselves to those $\kappa$ for which  $\kappa\leq |V_\zeta|$ for all $\zeta>\omega$, or equivalently to those $\kappa$ for which   $\kappa\leq |V_{\omega+1}|$. Given that $|V_{\omega+1}|=2^{\aleph_0}$, and the least possible value of the latter is $\aleph_1$, we conclude that the only acceptable values for $\kappa$ are $\aleph_0$ and $\aleph_1$, and therefore we should consider only  the classes
\begin{equation} \label{E:aleph0}
{\rm NT}_{\aleph_0}=\bigcup_{\zeta>\omega}{\rm NT}_{\aleph_0}(V_\zeta), \quad {\rm NT}_{\aleph_1}=\bigcup_{\zeta>\omega}{\rm NT}_{\aleph_1}(V_\zeta).
\end{equation}
So ${\rm NT}_{\aleph_0}$ and ${\rm NT}_{\aleph_1}$, with  ${\rm NT}_{\aleph_0}\subseteq {\rm NT}_{\aleph_1}$,  are the two natural candidate definitions for the class ${\rm NT}$ of nontypical sets. Both are expressible in the language of set theory, so they  satisfy condition (\ref{E:second}), and we shall see below that ${\rm OD}\subseteq {\rm NT}_{\aleph_0}\subseteq {\rm NT}_{\aleph_1}$, so in particular ${\rm Df}(V)\subseteq {\rm NT}_{\aleph_0}\subseteq {\rm NT}_{\aleph_1}$, i.e., they  satisfy also condition (\ref{E:first}).

Let us first deal with   ${\rm NT}_{\aleph_0}$. By Fact \ref{F:wideapp}, and using  Reflection, it follows that ${\rm NT}_{\aleph_0}$ is the class of algebraic elements of $V$, i.e., those belonging to finite sets which are definable in $V$ with ordinal parameters, i.e. to finite  OD sets.   Actually this class,  as well as the class ${\rm HNT}_{\aleph_0}=\{x:TC(\{x\})\subseteq {\rm NT}_{\aleph_0}\}$ of hereditarily nontypical sets, is not new. It has already come up and been investigated  by J.D. Hamkins and C. Leahy in \cite{HL16}, through a different motivation and   under the name class of ``ordinal algebraic'' sets,  denoted ${\rm OA}$. The only difference between  the definitions of ${\rm OA}$ and ${\rm NT}_{\aleph_0}$ is that while ${\rm NT}_{\aleph_0}=\bigcup_{\zeta>\omega}{\rm NT}_{\aleph_0}(V_\zeta)$, ${\rm OA}$ is defined in \cite{HL16} as the class of sets which are algebraic in the structures $V_\zeta$ by the extra help of ordinal parameters. This is equivalent to saying  that ${\rm OA}=\bigcup_{\zeta>\omega}{\rm NT}^*_{\aleph_0}(V_\zeta)$, where ${\rm NT}^*_{\aleph_0}(V_\zeta)$ denotes  the set  of algebraic elements of the structure $(V_\zeta,\in,\zeta)$, which is  $(V_\zeta,\in)$ endowed with ordinal parameters $<\zeta$.

Now it is well-known from  \cite{MS71}, that although for each  particular $V_\zeta$  definability in $(V_\zeta,\in,\zeta)$ and definability in $(V_\zeta,\in)$ do not coincide, the {\em totality} of sets definable in some $(V_\zeta,\in,\zeta)$ is no different than the totality of sets  definable in some $(V_\zeta,\in)$. This is due to the following  key fact (see the Extended Reflection Principle in \cite{MS71}).

\begin{Lem} \label{L:MS}
For any ordinals $\theta_1,\ldots,\theta_n$, there is an ordinal $\eta>\theta_1,\ldots,\theta_n$ such that $\theta_1,\ldots,\theta_n\in  {\rm Df}(V_\eta)$.
\end{Lem}
It is because of  this Lemma that  the class ${\rm OD}$ of ordinal definable sets can be defined just as $\bigcup_{\zeta\in On}{\rm Df}(V_\zeta)$, although we very often allow ordinal parameters in the definitions, i.e., we practically deal with the sets in $\bigcup_{\zeta\in On}{\rm Df}(V_\zeta,\zeta)$. For the same reason the following holds.

\begin{Lem} \label{L:similar}
${\rm OA}={\rm NT}_{\aleph_0}$.
\end{Lem}

{\em Proof.} We have to show that $\bigcup_{\zeta>\omega}{\rm NT}^*_{\aleph_0}(V_\zeta)=\bigcup_{\zeta>\omega}{\rm NT}_{\aleph_0}(V_\zeta)$. Since trivially for every $\zeta>\omega$, ${\rm NT}_{\aleph_0}(V_\zeta)\subseteq {\rm NT}^*_{\aleph_0}(V_\zeta)$, one inclusion is obvious.  For the converse, let  $a\in {\rm NT}^*_{\aleph_0}(V_\zeta)$, for some $\zeta$.  Then there are $\phi(x,\overline{y})$ and $\overline{\theta}<\zeta$  such that $a\in \phi(V_\zeta,\overline{\theta})$ and $|\phi(V_\zeta,\overline{\theta})|<\aleph_0$. By Lemma \ref{L:MS}, there is $\eta$ such that $\zeta,\overline{\theta}\in {\rm Df}(V_\eta)$. Then $V_\zeta\in V_\eta$ and by the absoluteness of satisfaction relation we have that for all $x\in V_\zeta$, $$V_\zeta\models\phi(x,\overline{\theta})\Leftrightarrow V_\eta\models (V_\zeta\models\phi(x,\overline{\theta})).$$ Since $\zeta,\overline{\theta}$ are definable in $V_\eta$, the formula  $\psi(x,\zeta,\phi,\overline{\theta}):=(V_\zeta\models\phi(x,\overline{\theta}))$ defines $\phi(V_\zeta,\overline{\theta})$ in $V_\eta$ without parameters, i.e., $\phi(V_\zeta,\overline{\theta})=\psi(V_\eta)$. Since $a\in \psi(V_\eta)$ and $|\psi(V_\eta)|=|\phi(V_\zeta,\overline{\theta})|<\aleph_0$, it follows that $a\in {\rm NT}_{\aleph_0}(V_\eta)$. This proves that $\bigcup_{\zeta>\omega}{\rm NT}^*_{\aleph_0}(V_\zeta)\subseteq \bigcup_{\zeta>\omega}{\rm NT}_{\aleph_0}(V_\zeta)$.  \telos

\vskip 0.2in

In view of Lemma \ref{L:MS} and its impact on the definition of ${\rm OD}$, the following simple characterizations of the classes ${\rm NT}_{\aleph_0}$ and ${\rm NT}_{\aleph_1}$ come out easily.

\begin{Lem} \label{L:chracter}
(i) A set  belongs to ${\rm NT}_{\aleph_0}$ iff it belongs to some finite  ${\rm OD}$ set.

(ii) A set  belongs to ${\rm NT}_{\aleph_1}$ iff it belongs to some countable ${\rm OD}$ set.
\end{Lem}

{\em Proof.} Let us sketch (ii). Let $a\in {\rm NT}_{\aleph_1}$. Then there are $\zeta$ and $\phi(x)$ such that  $a\in\phi(V_\zeta)$ and $|\phi(V_\zeta)|< \aleph_1$, i.e., $|\phi(V_\zeta)|\leq \aleph_0$. But clearly $\phi(V_\zeta)\in {\rm OD}$. Conversely let $A\in {\rm OD}$, $a\in A$ and $|A|\leq \aleph_0$. Then $A\in {\rm Df}(V_\zeta)$ for some $\zeta$. Since $A\subseteq V_\zeta$, clearly there is $\phi$ such that  $A=\phi(V_\zeta)$. Then $a\in \phi(V_\zeta)$ and  $|\phi(V_\zeta)|\leq \aleph_0$, therefore  $a\in {\rm NT}_{\aleph_1}(V_\zeta)\subseteq {\rm NT}_{\aleph_1}$. \telos

\vskip 0.2in

It follows immediately from (i) above that ${\rm OD}\subseteq {\rm OA}={\rm NT}_{\aleph_0}$.  However an interesting fact established in \cite{HL16} is that  the class ${\rm HOA}$ of hereditarily ordinal algebraic sets coincides with the class ${\rm HOD}$ of hereditarily ordinal definable  sets. Let  ${\rm HNT}_{\aleph_i}$ denote the hereditary subclasses of ${\rm NT}_{\aleph_i}$, for $i=0,1$ respectively. By \ref{L:similar}, ${\rm HOA}={\rm HNT}_{\aleph_0}$.

\begin{Thm} \label{T:HL}
{\rm (\cite{HL16})} ${\rm HOA}={\rm HOD}$. Therefore also ${\rm HNT}_{\aleph_0}={\rm HOD}$.
\end{Thm}
(On the other hand it is known that   ${\rm NT}_{\aleph_0}\neq {\rm OD}$ is consistent. Specifically it was proved in \cite{GL87} that there is a generic extension of $\textbf{L}$  containing  a pair  $\{X,Y\}$ which is ${\rm OD}$, while  neither $X$ nor $Y$ is ${\rm OD}$. Thus $X,Y\in {\rm NT}_{\aleph_0}\backslash {\rm OD}$.)

By the preceding theorem,   ${\rm NT}_{\aleph_0}$ loses some of its interest as a   candidate class for the definition of ${\rm NT}$, since its hereditary subclass collapses to the familiar inner model ${\rm HOD}$. So if one is looking  for a really {\em new} inner model of ${\rm ZF}$ which  strictly exceeds  ${\rm HOD}$, I think the only option left  is to  identify  ${\rm NT}$ with  the class ${\rm NT}_{\aleph_1}$. It turns out that the subclass ${\rm HNT}_{\aleph_1}$ of the latter is indeed a new inner model of ${\rm ZF}$.

\begin{Thm} \label{T:simplified}
${\rm ZFC}$ proves that  ${\rm HNT}_{\aleph_1}$ is an inner model of ${\rm ZF}$ such that ${\rm HOD}\subseteq {\rm HNT}_{\aleph_1}$.
\end{Thm}

{\em Proof.} We work in ${\rm ZFC}$. That ${\rm HOD}\subseteq {\rm HNT}_{\aleph_1}$ follows from the fact that already ${\rm HOD}\subseteq {\rm HNT}_{\aleph_0}$ (actually ${\rm HOD}={\rm HNT}_{\aleph_0}$ by \ref{T:HL}) and ${\rm HNT}_{\aleph_0}\subseteq{\rm HNT}_{\aleph_1}$.

Let us write for simplicity ${\rm HNT}$ throughout this proof instead of  ${\rm HNT}_{\aleph_1}$.  Extensionality holds in ${\rm HNT}$ because of the transitivity of the latter, and Foundation is true  trivially because is true  in the underlying universe $V$. Also Infinity holds trivially since $\omega\in{\rm HOD}\subseteq {\rm HNT}$. So it remains to  prove Pairing, Union, Powerset and Replacement. The proof is based on the  characterization given in Lemma \ref{L:chracter} (ii), that $a\in {\rm HNT}={\rm HNT}_{\aleph_1}$ if and only if there is $A\in {\rm OD}$ such that $a\in A$ and $|A|\leq \aleph_0$.

\vskip 0.1in
{\em Pairing.} Let $a,b\in {\rm HNT}$. We have to show that $\{a,b\}\in {\rm NT}$. By assumption there are $A,B\in {\rm OD}$ such that $a\in A$, $b\in B$ and $|A|,|B|\leq \aleph_0$. Let $C=\{\{x,y\}:x\in A, y\in B\}$. Clearly $C\in {\rm OD}$,  $|C|\leq \aleph_0$ and $\{a,b\}\in C$.

\vskip 0.1in

{\em Union.} Let $a\in {\rm HNT}$. It suffices to see that $\cup a\in {\rm NT}$. Let $a\in A$, where $A\in {\rm OD}$ such that  $|A|\leq \aleph_0$. Let $B=\{\cup x:x\in A\}$. Clearly $B\in {\rm OD}$, $|B|\leq \aleph_0$ and $\cup a\in B$.

\vskip 0.1in

{\em Powerset.} Let $a\in {\rm HNT}$. It suffices to show that ${\cal P}^{\rm HNT}(a)={\cal P}(a)\cap {\rm HNT}$ belongs to NT. Let   $a\in A$ for some $A\in {\rm OD}$ with $|A|\leq \aleph_0$.   For every $x$,  ${\cal P}(x)\cap {\rm HNT}$ is a set (by Separation in $V$), and (by Replacement)  so is also  $B=\{{\cal P}(x)\cap {\rm HNT}:x\in A\}$. Moreover it is easy to check that $B\in {\rm OD}$. Since obviously
$|B|\leq |A|\leq \aleph_0$ and ${\cal P}(a)\cap {\rm HNT}\in B$, we are done.

\vskip 0.1in

{\em Replacement.} Let $a\in {\rm HNT}$ and let $\phi(x,y,b_1,\ldots,b_n)$ be a formula with parameters  $b_i\in {\rm HNT}$, such that ${\rm HNT}\models (\forall x\in a)(\exists !y)\phi(x,y,b_1,\ldots,b_n)$, or equivalently, $(\forall x\in a)(\exists !y)\phi^{\rm HNT}(x,y,b_1,\ldots,b_n)$, where $\phi^{\rm HNT}$ is the usual relativization of $\phi$ to the class ${\rm HNT}$.  $\phi^{\rm HNT}$ defines a functional relation on $a$, so let us write $F^{\rm HNT}_{\phi(\overline{b})}(x)=y$  instead of $\phi^{\rm HNT}(x,y,b_1,\ldots,b_n)$, where $\overline{b}=\langle b_1,\ldots,b_n\rangle$. If for some tuple $\langle b_1,\ldots,b_n\rangle$, $\phi(x,y,b_1,\ldots,b_n)$ does not define a function, we set $F^{\rm HNT}_{\phi(\overline{b})}(x)=\emptyset$.

Under this notation we have to show that for the given $a,\overline{b}\in {\rm HNT}$, the set $c=F^{\rm HNT}_{\phi(\overline{b})}[a]=\{F^{\rm HNT}_{\phi(\overline{b})}(x):x\in a\}$ is an element of ${\rm HNT}$. By our assumption there are  countable $A,B_1,\ldots,B_n\in{\rm OD}$ such that  $a\in A$ and  $b_i\in B_i$, for $i=1,\ldots,n$. Since  $F^{\rm HNT}_{\phi(\overline{b})}$ is a function within ${\rm HNT}$, for each $x\in a$, $F^{\rm HNT}_{\phi(\overline{b})}(x)\in {\rm HNT}$, therefore $c\subseteq {\rm HNT}$. So it suffices to show that $c\in {\rm NT}$, i.e., $c\in C$ for some countable  $C\in {\rm OD}$. Let
$$C=\{F^{\rm HNT}_{\phi(\overline{w})}[z]:z\in A, w_1\in B_1,\ldots,w_n\in B_n\}.$$
Using Reflection and the fact that $A,B_i$ are in ${\rm OD}$, it is not hard to see that $C\in {\rm OD}$. Since $a\in A$ and $b_i\in B_i$, we have that $c=F^{\rm HNT}_{\phi(\overline{b})}[a]$ belongs to $C$. Moreover, since the variables $z$ and $w_i$ range over the countable sets $A$ and $B_i$, respectively, it follows that   $$|C|\leq |A\times B_1\times\cdots\times B_n|\leq\aleph_0.$$
Thus $C$ is a countable  ${\rm OD}$ set containing $c$ and we are done. \telos

\vskip 0.2in

We do not know if $AC$ holds in ${\rm HNT}_{\aleph_1}$ or not. We only know that  we cannot prove in ${\rm ZFC}$ that $AC$ fails in ${\rm HNT}_{\aleph_1}$. Because if ${\rm ZFC}$ is consistent, then so is  ${\rm ZFC}+V={\rm HOD}$. But the latter theory implies  ${\rm HOD}={\rm HNT}_{\aleph_1}$ and ${\rm HOD}$ satisfies $AC$, so ${\rm ZFC}+AC^{{\rm HNT}_{\aleph_1}}$ is consistent.  The latter would be also  a consequence of the consistency of ${\rm ZFC}+V={\rm HNT}_{\aleph_1}$ alone,  no matter whether ${\rm HOD}={\rm HNT}_{\aleph_1}$ or ${\rm HOD}\neq {\rm HNT}_{\aleph_1}$.

In general, in view of the inclusions ${\rm HOD}\subseteq {\rm HNT}_{\aleph_1}\subseteq V$, the questions  regarding the  consistency of the various mutual relationships among these three classes arise naturally. The simplest such relationship is  of course ${\rm HOD}={\rm HNT}_{\aleph_1}=V$ and  follows from  $V={\rm HOD}$, whose consistency is well-known.\footnote{Note by the way that, as is the case with the classes ${\rm HOD}$ and ${\rm OD}$ (described in \cite[p. 276]{MS71}), the following equivalences hold: ${\rm HNT}_{\aleph_1}={\rm NT}_{\aleph_1}$$\Leftrightarrow$$V={\rm HNT}_{\aleph_1}$$\Leftrightarrow$$V={\rm NT}_{\aleph_1}$. The last equivalence,  as well as  $\Leftarrow$ of the first equivalence are obvious. Concerning $\Rightarrow$ of the first equivalence, assume ${\rm HNT}_{\aleph_1}={\rm NT}_{\aleph_1}$. For every $\alpha\in Ord$, clearly  $V_\alpha\in {\rm OD}\subseteq {\rm NT}_{\aleph_1}$, so $V_\alpha\in {\rm HNT}_{\aleph_1}$, whence  $V_\alpha\subseteq {\rm HNT}_{\aleph_1}$ and therefore $V={\rm HNT}_{\aleph_1}$.}  Concerning the other ones we have only two  partial answers.

Perhaps the most  urgent question  to answer is the  {\em existence} itself of typical sets, i.e., the consistency of $V\neq {\rm NT}_{\aleph_1}$. For if we do not
know whether $V$ can be separated from ${\rm NT}_{\aleph_1}$, the definition of ${\rm TP}_{\aleph_1}$ is vacuous.
Fortunately this question can be affirmatively settled by a lot of forcing notions which have  rich sets of automorphisms.  More specifically  the following holds.

\begin{Thm} \label{T:affirm}
(i) Let $M\models{\rm ZFC}$,  $\P\in M$ be a forcing notion, and $G\subseteq \P$ be $M$-generic. For any  $p\in \P$, let $A_p={\rm Aut}^M_{\{p\}}(\P)$ be the set of automorphisms of $\P$ in $M$ which fix $p$. Assume further that for  every $p\in G$,  the set $\{\pi''G:\pi\in A_p\}$ is uncountable in $M[G]$. Then   $G$ is  typical in $M[G]$, and hence $M[G]\models V\neq {\rm NT}_{\aleph_1}$.

(ii) In particular, there are $M$ and $G$ such that $M[G]\models |{\cal P}(\omega)\cap {\rm TP}_{\aleph_1}|=2^{\aleph_0}$ (where ${\rm TP}_{\aleph_1}=V\backslash {\rm NT}_{\aleph_1}$), i.e. $M[G]$ contains continuum many typical reals.
\end{Thm}

{\em Proof.} (i) First note that the set $A_p$ belongs to $M$, and hence to $M[G]$, so $\{\pi''G:\pi\in A_p\}$ is also an element of $M[G]$. Next it is known that  there is an abundance of forcing notions $\P$ satisfying the requirement of the theorem. For example such is the poset  $\P$ of finite functions $p$ with $dom(p)\subset \omega$ and $rng(p)\subseteq \{0,1\}$, ordered by reverse inclusion, which adds a single Cohen real.  The automorphisms of $\P$ are induced by the permutations $\pi:\omega\rightarrow \omega$, and are defined as follows: for every $p\in\P$, $dom(\pi(p))=\pi[dom(p)]$ and  $\pi(p)(\pi(n))=p(n)$. Then  for every generic $G$ and $p\in G$, $\{\pi''G:\pi\in A_p\}$ is uncountable in $M[G]$.

To verify the claim of the theorem, let $\P\in M$ and $G$ be as  stated.  We have to show, according to Lemma \ref{L:chracter} (ii), that for every $A\in {\rm OD}^{M[G]}$ such that $G\in A$,  $M[G]\models |A|>\aleph_0$. Pick any $A\in {\rm OD}^{M[G]}$ containing $G$. Then there are a formula $\phi(x,y_1,\ldots,y_n)$ and ordinals $\alpha_1,\ldots,\alpha_n$ such that $A=\{x:M[G]\models\phi(x,\alpha_1,\ldots,\alpha_n)\}$ and  $M[G]\models\phi(G,\alpha_1,\ldots,\alpha_n)$. The key fact here is that  all generic subsets $G$  of $\P$ have a ``common'' $\P$-name, sometimes called ``canonical name'', namely $\Gamma=\{\langle \breve{q},q\rangle:q\in \P\}$. If $t^G$ denotes the $G$-interpretation of a $\P$-name $t$ into $M[G]$, then for every generic $G$, $\Gamma^G=G$. Now $M[G]\models\phi(G,\alpha_1,\ldots,\alpha_n)$ means that there is some  $p\in G$ such that
$$p\parallel\!\!\!-\phi(\Gamma,\breve{\alpha}_1,\ldots,\breve{\alpha}_n),$$ where $\Gamma$ is the canonical name.  Fix such a $p\in G$. It is  well-known that for every automorphism  $\pi$ and generic $G$, $\pi''G$ is  generic too and moreover $M[G]=M[\pi''G]$. But then for every  $\pi\in A_p$, $p\in \pi''G$, so $p\parallel\!\!\!\!-\phi(\Gamma,\breve{\alpha}_1,\ldots,\breve{\alpha}_n)$ implies also that $M[\pi''G]\models\phi(\Gamma^{\pi''G},\alpha_1,\ldots,\alpha_n)$, or  $M[\pi''G]\models\phi(\pi''G,\alpha_1,\ldots,\alpha_n)$. Therefore for every $\pi\in A_p$,
$$M[G]=M[\pi''G]\models\phi(\pi''G,\alpha_1,\ldots,\alpha_n).$$
It follows that the set $A$, which is the extension of $\phi$ in $M[G]$, has as subset the set   $\{\pi''G:\pi\in A_p\}$ which is uncountable in $M[G]$ by assumption, and thus also $M[G]\models|A|>\aleph_0$.

(ii) Just take $\P$ to be the forcing notion mentioned in the beginning of the proof for adding a Cohen real, so, essentially,  $G\in {\cal P}(\omega)$. Then $M\models|Aut(\P)|=2^\omega$ and it is a folklore fact that there are continuum many images $\pi(G)$ of $G$ in  $M[G]$ all of which are generic sets, hence typical. So   $M[G]\models |{\cal P}(\omega)\cap {\rm TP}_{\aleph_1}|=2^{\aleph_0}$. \telos

\vskip 0.2in

One of the referees informed me that the fact shown in  Theorem \ref{T:affirm}, i.e., that a generic set which satisfies the given conditions does not belong to ${\rm HNT}_{\aleph_1}$,  has  been established also in \cite{KL18} for Cohen and Solovay-random extensions by a different and more complex argument.

Nevertheless, not all generic sets are typical. This is a side consequence of the following relevant result of V.G. Kanovei and V.A. Lyubetsky in \cite{KL17}, which implies that the classes ${\rm HOD}$ and ${\rm HNT}_{\aleph_1}$ can be separated.

\begin{Thm} \label{T:KL}
{\rm (\cite[Theorem 4]{KL17})} There is a generic extension ${\rm \mathbf{L}}[(x_n)_{n<\omega}]$ of the constructible universe ${\rm \mathbf{L}}$ by a sequence of reals $x_n\in 2^\omega$, in which it is true that $\{x_n:n<\omega\}$ is a countable $\Pi^1_2$ set with no ${\rm OD}$ elements.
\end{Thm}
In the proof of this  theorem the sets $x_n$ are added to ${\rm \mathbf{L}}$ {\em generically}, so, in contrast to Theorem \ref{T:affirm}, generics are used here, essentially, to show the existence of {\em nontypical} sets which are not ordinal definable. Namely, the following holds.

\begin{Cor} \label{C:follows}
If ${\rm ZF}$ is consistent, then so is ${\rm ZFC}+{\rm HOD}\neq {\rm HNT}_{\aleph_1}$. In particular it is consistent that ${\cal P}(\omega)\cap {\rm OD}\varsubsetneq{\cal P}(\omega)\cap {\rm NT}_{\aleph_1}$.
\end{Cor}

{\em Proof.} If $X=\{x_n:n<\omega\}$ is the set of Theorem \ref{T:KL}, clearly $X\in {\rm OD}$ and $X$ is countable. Therefore every $x_n$ belongs to ${\rm NT}_{\aleph_1}$. Moreover every $x_n $ is a real, so $x_n\in {\rm HNT}_{\aleph_1}$.  On the other hand, since $x_n\notin {\rm OD}$, it follows that  $x_n\in {\rm HNT}_{\aleph_1}\backslash {\rm HOD}$. By the same token,  $x_n\in ({\cal P}(\omega)\cap {\rm NT}_{\aleph_1})\backslash ({\cal P}(\omega)\cap{\rm OD})$. \telos

\begin{Que} \label{Q:two}
Are the following  consistent with  ${\rm ZFC}$?

1) ${\rm HOD}\neq {\rm HNT}_{\aleph_1}\neq V$,

2) ${\rm HOD}={\rm HNT}_{\aleph_1}\neq V$,

3) ${\rm HOD}\neq {\rm HNT}_{\aleph_1}= V$.
\end{Que}

\begin{Que} \label{Q:three}Is it consistent with ${\rm ZFC}$ that $AC$ fails in ${\rm HNT}_{\aleph_1}$?
\end{Que}

\vskip 0.1in

\textbf{Acknowledgements} I am  indebted to two anonymous referees for some serious    corrections  and   several  clarifications and suggestions which significantly  improved this article.

\end{document}